\documentclass[reqno]{amsart} \usepackage{amscd}
\usepackage{epsf}
\newtheorem{theorem}{Theorem}[section]
\newtheorem{proposition}[theorem]{Proposition}
\newtheorem{lemma}[theorem]{Lemma}

\theoremstyle{remark} 
\numberwithin{equation}{section}
\newcommand{\field}[1]{\ensuremath{\mathbb{#1}}}

\newcommand{\RR}{\field{R}}

\renewcommand{\d}{\operatorname{d}}

\newcommand{\del}{\partial}

\begin{document}
\title[A Variational Proof For  The Existence]
{A variational proof for the existence of a conformal metric with preassigned
negative Gaussian curvature for compact Riemann surfaces of genus $>1$}
\author{Rukmini Dey} \address{School of Mathematics\\
Harish-Chandra-Research Institute, Jhusi, Allahabad, India}
\email{rkmn@mri.ernet.in}
\begin{abstract}
Given an smooth function $K <0$ we prove a result by Berger, Kazhdan and
others that in every conformal class there exists a metric which attains this
function as its Gaussian curvature for a compact Riemann surface of  genus
$g>1$. We do so by minimizing an appropriate functional using elementary
analysis. In particular for $K$ a negative constant, this provides an
elementary proof of the  uniformization theorem for compact Riemann surfaces
of genus $g >1$.
\end{abstract}
\maketitle

\section{Introduction}

In this paper we present a variational proof of a result by Berger ~\cite{B},
Kazhdan , Warner ~\cite{KW} and others, namely given  an arbitrary smooth
function $K <0$ we show that in every conformal class there exists a metric
which attains this function as its Gaussian curvature for a compact Riemann
surface of  genus $g>1$. In particular, this result includes the
uniformization theorem of H. Poincar\'{e} ~\cite{P} when $K$ is a negative
 constant. In his proof Berger  considers the critical points of  a
 functional subject to the  Gauss-Bonnet condition. He  shows that the
functional is bounded from below and uses the Friedrich's
inequality to complete the proof. The functional we choose is
positive definite so that it is automatically bounded from below.
Our proof is elementary, using Hodge theory, i.e., the existence
of the Green's operator for the Laplacian. Our proof could be
useful for analysing the appropriate condition on $K$ for a
corresponding result for genus $g=1$ and $g=0$ ~\cite{KW},
~\cite{XY}, ~\cite{CL}, the two other cases considered by Berger,
Kazdan and Warner. Another variational proof of the uniformization
theorem for genus $g>1$ can  be found in a gauge-theoretic context
in ~\cite{H} which uses Uhlenbeck's weak compactness theorem for
connections with $L^p$ bounds on curvature ~\cite{U}.

Let $M$ be a compact Riemann surface of
genus $g>1$ and let $\d s^2=h\d z \otimes\d\bar{z}$ be a  metric
on $M$ normalized such that the total area of $M$ is $1$. Let $K<0$.
We  minimize the functional
\begin{equation*}
S(\sigma)=\int_M (K(\sigma)-K)^2 e^{2\sigma}d\mu
\end{equation*}
over $W^{2,2}(M)$, where $K(\sigma)$ stands for the Gaussian curvature
of the  metric $e^{\sigma}\d s^2$, and $d\mu = \frac{\sqrt{-1}}{2} h \d z
\wedge\d\bar{z}$ is the area form for the metric $ds^2$. Using Sobolev
embedding theorem we show that
$S(\sigma)$ takes its absolute minimum on $W^{2,2}(M)$ which corresponds to
a $C^{\infty}$ metric on $M$ of negative curvature $K$.

\section{The main theorem}
\subsection{}
All notations are as in section 1.

The functional $ S(\sigma)=\int_M(K(\sigma)-K)^2e^{2\sigma} d \mu $
is non-negative on $W^{2,2}(M)$, so that its
infimum
\begin{equation*}
S_0 = \inf\{S(\sigma),\,\sigma \in W^{2,2}(M)\}
\end{equation*}
exists and is non-negative. Let
$\{\sigma_n\}_{n=1}^{\infty}\subset W^{2,2}(M)$ be a
corresponding minimizing sequence,
\begin{equation*}
\lim_{n\rightarrow\infty} S(\sigma_n) =S_0.
\end{equation*}
Our main result is the following
\begin{theorem} Let $M$ be a compact Riemann surface of genus $g >1$.
The infimum $S_0$ is attained at $\sigma\in C^{\infty}(M,\RR)$,
i.e.~the minimizing sequence $\{\sigma_n\}$ contains a subsequence that
converges in $W^{2,2}(M) $ to $\sigma\in C^{\infty} (M,\RR)$ and
$S(\sigma)=0$. The corresponding metric $e^{\sigma}h\d z\otimes\d\bar{z}$ is
the unique metric on $M$ of negative curvature $K$.
\end{theorem}

\subsection{Uniform bounds}
 Since $\{\sigma_n\}$ is a minimizing sequence, we have the obvious
inequality
\begin{equation}\label{1}
 S[\sigma_n]=\int_{M}(K_n-K)^2e^{2\sigma_n}d\mu
=\int_M (K_0 - \frac{1}{2}\Delta_h\sigma_n -Ke^{\sigma_n})^2 d\mu\leq m
\end{equation}
for some $m>0$, where we denoted by $K_n$ the Gaussian curvature $K(\sigma_n)$
of the metric $e^{\sigma_n}h$ and by $K_0$ that of the metric $h$, and used
that
\begin{equation*}
K_n=e^{-\sigma_n}(K_0 - \frac{1}{2} \Delta_h \sigma_n).
\end{equation*}

Note: Here $\Delta_h = 4h^{-1}\frac{\del^2}{\del z\del\bar{z}} $
stands for the Laplacian defined by the metric $h$ on $M$.
\begin{lemma} There exist a constant $C_1$  such that, uniformly in
$n$,
\begin{eqnarray*}
 \int_M (\Delta_h\sigma_n)^2 d\mu <C_1
\end{eqnarray*}
\end{lemma}

\begin{proof} By Minkowski inequality, and using \eqref{1}, we get
\begin{eqnarray*}
[\int_M (-\frac{1}{2} \Delta_h\sigma_n -Ke^{\sigma_n})^2 d \mu]^{1/2} &\leq&
[\int_M (K_0 - \frac{1}{2}\Delta_h\sigma_n -Ke^{\sigma_n})^2 d \mu]^{1/2}\\
&+& [\int_M (K_0)^2 d \mu ]^{1/2}\leq m^{1/2}+c=C,
\end{eqnarray*}
so that
\begin{equation}\label{2}
\frac{1}{4}\int_M (\Delta_h\sigma_n)^2 d\mu+ \int_M K^2 e^{2\sigma_n}d\mu +
\int_M \Delta_h\sigma_n e^{\sigma_n} K d \mu \leq C^2.
\end{equation}

We will show that
\begin{equation}
\int_M K^2 e^{2\sigma_n}d\mu + \int_M \Delta_h\sigma_n e^{\sigma_n} K d \mu
= B_{n1} + B_{n2} + B_{n3}
\end{equation}
where $B_{n1} \geq 0$, $B_{n2} \geq 0$, $|B_{n3}| \leq 3 D^2$ where $D^2$ is
a constant independent of $n$.

From (2.3) and (2.2) the result  will follow 
since we will have
\begin{eqnarray*}
C^2 + 3D^2 &\geq& C^2 - B_{n3} \geq \frac{1}{4}\int_M (\Delta_h\sigma_n)^2 d\mu
+ B_{n1} + B_{n2} \\
&\geq& \frac{1}{4}\int_M (\Delta_h\sigma_n)^2 d\mu
\end{eqnarray*}
Just renaming the constants, we will have the result.

Integrating by parts we get
\begin{eqnarray*}
 \int_M \Delta_h\sigma_n e^{\sigma_n} K d \mu &=& - \int_M |\partial_{z} \sigma_n|^2 e^{\sigma_n} K d \mu -  \int_M (\partial_z \sigma_n) (\partial_{\bar{z}} K) e^{\sigma_n} d \mu \\
& & -\int_{\partial M} (\partial_{\nu} \sigma_n) K e^{\sigma_n} d \mu \\
&=&  \int_M |\partial_{z} \sigma_n|^2 e^{\sigma_n} |K|d \mu -  \int_M (\partial_z \sigma_n) g |K| e^{\sigma_n} d \mu  
\end{eqnarray*}
since $K$ is negative, $\partial_{\nu} \sigma_{n} |_{\partial M} = 0$ and 
where we define $g= \frac{\partial_{\bar{z}} K}{|K|}$.

Let $M = \Omega_{n1} \cup \Omega_{n2} \cup \Omega_{n3}$, a disjoint union of 
sets 
defined as follows:

On $\Omega_{n1}$,  $(1)$ $ |\partial_{z} \sigma_n |> |g|.$

On $\Omega_{n2}$, $(2)$ $ |\partial_{z} \sigma_n |\leq |g|$ and
$ |K| e^{\sigma_n} > |g|^2.$

On $\Omega_{n3}$, $(3)$  $ |\partial_{z} \sigma_n |\leq |g|$ and
$ |K| e^{\sigma_n} \leq |g|^2.$

Let $B_{ni} = \int_{\Omega_{ni}} K^2 e^{2\sigma_n}d\mu + \int_{\Omega_{ni}} 
\Delta_h\sigma_n e^{\sigma_n} K d \mu$, $i=1,2,3$.

We will show that $B_{n1} \geq 0$.
\begin{eqnarray*} 
B_{n1} &=& \int_{\Omega_{n1}} K^2 e^{2\sigma_n} d \mu + \int_{\Omega_{n1}} |\partial_{z} \sigma_n|^2 |K| e^{\sigma_n} d \mu 
- \int_{\Omega_{n1}} (\partial_{z} \sigma_n) g |K| e^{\sigma_n} d \mu\\
&\geq&  \int_{\Omega_{n1}} K^2 e^{2\sigma_n} d \mu + \int_{\Omega_{n1}} |\partial_{z} \sigma_n|^2 |K| e^{\sigma_n} d \mu 
- \int_{\Omega_{n1}} |\partial_{z} \sigma_n| |g| |K| e^{\sigma_n} d \mu\\
&=&\int_{\Omega_{n1}} K^2 e^{2\sigma_n} d \mu + \int_{\Omega_{n1}} |\partial_{z} \sigma_n| |K| e^{\sigma_n}(|\partial_z \sigma_n| - |g|) d\mu \\
&\geq& 0 
\end{eqnarray*}
by $(1)$ in the definition of $\Omega_{n1}$.

Next we shall show that $B_{n2} \geq 0$.
\begin{eqnarray*}
B_{n2} &=& \int_{\Omega_{n2}} K^2 e^{2\sigma_n} d \mu + \int_{\Omega_{n2}} |\partial_{z} \sigma_n|^2 |K| e^{\sigma_n} d \mu - \int_{\Omega_{n2}} (\partial_{z} \sigma_n) g |K| e^{\sigma_n} d \mu\\
&\geq&  \int_{\Omega_{n2}} K^2 e^{2\sigma_n} d \mu + \int_{\Omega_{n2}} |\partial_{z} \sigma_n|^2 |K| e^{\sigma_n} d \mu - \int_{\Omega_{n2}} |\partial_{z} \sigma_n| |g| |K| e^{\sigma_n} d \mu\\
&=& \int_{\Omega_{n2}} |K| e^{\sigma_n} (|K|e^{\sigma_n} - |\partial_{z} \sigma_n| |g|) d \mu + \int_{\Omega_{n2}} |\partial_z \sigma_n|^2 e^{\sigma_n} |K| d \mu \\
&\geq& \int_{\Omega_{n2}} |K| e^{\sigma_n} (|K|e^{\sigma_n} - |g|^2) d \mu + 
\int_{\Omega_{n2}} |\partial_z \sigma_n|^2 e^{\sigma_n} |K| d \mu \\
&\geq& 0,
\end{eqnarray*}
by using the two conditions $(2)$ defining $\Omega_{n2}$.

Next we shall show that $B_{n3}$ is uniformly bounded.
\begin{eqnarray*}
|B_{n3}| &\leq& \int_{\Omega_{n3}} K^2 e^{2\sigma_n} d \mu + \int_{\Omega_{n3}} |\partial_{z} \sigma_n|^2 |K| e^{\sigma_n} d \mu + \int_{\Omega_{n3}} |\partial_{z} \sigma_n| |g| |K| e^{\sigma_n} d \mu\\
&\leq& 3 D^2,
\end{eqnarray*}
where $D^2 = max |g|^4 \mu(M) $, where $\mu(M)$ is the volume of $M$. 
This  follows from the two conditions 
  $|K| e^{\sigma_n} \leq |g|^2$
and $ |\partial_{z} \sigma_n| \leq |g|$ on $\Omega_{n3}$.
$D^2$ is a finite constant ($max |g|^4$ is finite since  $K$  is
non-zero and the volume of $M$ is finite), independent of $n$.
Thus the result follows.
\end{proof}

\subsection{Pointwise convergence of zero mean-value part}

Next, for $\sigma\in C^{\infty}(M)$ denote by $m(\sigma)$ its mean value,
\begin{equation*}
m(\sigma)=\int_M\sigma d\mu,
\end{equation*}
and by $\tilde{\sigma}=\sigma-m(\sigma)$ denote its zero-mean value part.
For the minimizing sequence $\{\sigma_n\}$ we denote the corresponding
mean values by $m_n$.
(Note: we had normalised the volume $\int_M d \mu =1$.)

\begin{lemma}
The mean-value-zero part $\{\tilde{\sigma}_n\}_{n=1}^{\infty}$ of the
minimizing sequence

$\{\sigma_n\}_{n=1}^{\infty}$ is uniformly bounded in the Sobolev space $W^{2,2}(M)$.
\end{lemma}
\begin{proof}
By Hodge theory, there exists an operator $G$ such that
$G\Delta_h=I-P$, where $I$ is the identity operator in $L^2(M)$
and $P$ is the orthogonal projection onto kernel of $\Delta_h$. We
also know $\Delta_h : W^{2,2} \rightarrow L^2 $ boundedly and
$G:L^2 \rightarrow W^{2,2}$ is a bounded operator.

Now, by lemma 2.2 $ \{ \Delta_h \sigma_n \}$ are bounded uniformly
in $L^2$.

Thus, $ \{ G \Delta_h \sigma_n \}$ are bounded uniformly in
$W^{2,2}$.

But $ G \Delta_h \sigma_n = (I - P) \sigma_n = \tilde{\sigma_n} $.

\end{proof}
Now we can formulate the main result of this subsection.

\begin{proposition} The sequence $\{\tilde{\sigma}_n\}_{n=1}^{\infty}$
contains a subsequence $\{\tilde{\sigma}_{l_n}\}_{n=1}^{\infty}$ with
the following properties.
\begin{itemize}
\item[(a)] The sequences $\{\tilde{\sigma}_{l_n}\}_{n=1}^{\infty}$ and
$\{e^{\tilde{\sigma}_{l_n}+m_{l_n}}\}$ converge in $W^{2,2}(M)$
to continuous functions $\tilde{\sigma}$ and $u$
respectively. Moreover, $\tilde{\sigma}\in W^{2,2}(M)$.

\item[(b)]  The subsequence $\{\Delta_h \tilde{\sigma}_{l_n} \}$ converges
weakly in $L^2$ to $f\doteq \Delta_h^{distr} \tilde{\sigma}$---a
distribution Laplacian of $\tilde{\sigma}$.
\item[(c)] Passing to this subsequence $\{\tilde{\sigma}_{l_n} \}$, the
following limits exist
\begin{eqnarray*}
& \lim_{n\rightarrow\infty}\|\Delta_h \tilde{\sigma_{l_n}}\|_2=\|
\Delta_h^{distr} \tilde{\sigma}\|_2, \\ &
\lim_{n\rightarrow\infty}S(\sigma_{l_n})=S_0=\int_M(K_0-\frac{1}{2}
\Delta_h^{distr}\tilde{\sigma} -Ku)^2 d\mu.
\end{eqnarray*}
where
\begin{equation*}
\lim_{n\rightarrow\infty}e^{\tilde{\sigma}_{l_n}+ m_n}=u.
\end{equation*}
Infact, the convergence in (b) is strong in $L^2$.
\end{itemize}
\end{proposition}

\begin{proof} Part (a) follows from the Sobolev embedding theorem and Rellich
lemma since, for $\dim M = 2$, the space $W^{2,2}(M)$ is compactly
embedded into $C^0(M)$ (see, e.g.~\cite{A,K}). Therefore the
sequence $\{\tilde{\sigma}_n\}$  , which , according to lemma 2.3,
is uniformly bounded in $W^{2,2}(M)$, contains a convergent
subsequence in $C^0(M)$. Passing to this subsequence $\{
\tilde{\sigma}_{l_n} \}$ we can assume that there exists
mean-value zero function $\tilde{\sigma}\in C^0(M)$ such that
\begin{equation*}
\lim_{n\rightarrow\infty}\tilde{\sigma}_{l_n}=\tilde{\sigma}.
\end{equation*}
Since $\tilde{\sigma}_n$'s are uniformly bounded in a Hilbert space $W^{2,2}
(M)$, they weakly converge to $s\in W^{2,2}(M)$ (after passing to a
subsequence if necessary). The uniform limit coincides with $s$  so that
$\tilde{\sigma}=s\in W^{2,2}(M)$.

We have to show that $m_n=\int_M\sigma_n d\mu$ is bounded from above.
Suppose not, i.e. $m_{n} \rightarrow \infty$.
From the proof of  lemma (2.2)  since $\int_M (\Delta \sigma_n)^2 d \mu$ is 
positive $\int_M K^2 e^{2\sigma_n}d\mu + \int_M \Delta_h\sigma_n e^{\sigma_n} K d \mu \leq C^2 $, uniformly in $n$.

\begin{eqnarray*}
\int_M K^2 e^{2\sigma_n}d\mu &+& \int_M \Delta_h\sigma_n e^{\sigma_n} K d \mu \\
&=& e^{2 m_{n}} [\int_M K^2 e^{2 \tilde{\sigma}_n} d \mu + e^{-m_{n}} \int_M \Delta_h \tilde{\sigma}_n e^{\tilde{\sigma}_n} K d\mu] \leq C^2
\end{eqnarray*}

Now $e^{m_n} \rightarrow \infty$ as $n \rightarrow \infty$. But 
the right hand side of the previous equality is bounded from above, therefore
$\int_M K^2 e^{2 \tilde{\sigma}_n} d \mu + e^{-m_{n}} \int_M \Delta_h \tilde{\sigma}_n e^{\tilde{\sigma}_n} K d\mu$ tends to $0$.

Note that $\Delta_h \sigma_n = \Delta_h \tilde{\sigma}_n$. 
Let us abbreviate $A_n = \int_M K^2 e^{2 \tilde{\sigma}_n}d \mu$
\begin{eqnarray*}
|\int_M \Delta_h\tilde{\sigma}_n e^{\tilde{\sigma}_n} K d \mu| \leq  (\int_M (\Delta_h \tilde{\sigma}_n)^2 d \mu)^{\frac{1}{2}} A_n^{\frac{1}{2}} \leq C_1^{\frac{1}{2}} A_n^{\frac{1}{2}}
\end{eqnarray*}
where by lemma (2.2), $(\int_M (\Delta_h\tilde{\sigma}_n)^2 d \mu)^{\frac{1}{2}} \leq C_1^{\frac{1}{2}}$, where $C_1$ is  independent of $n$.
\begin{eqnarray*}
C^2 &\geq& e^{2 m_n} [ A_n + e^{-m_n} \int_M \Delta_h \tilde{\sigma}_n K e^{\tilde{\sigma}_n}d\mu ]\\
&\geq& e^{2m_n}[A_n - e^{-m_n} C_1^{\frac{1}{2}} A_n^{1/2}] \\
&=&e^{2m_n} (A_n)^{1/2}[ A_n^{1/2}  - e^{-m_n} C_1^{\frac{1}{2}}]\\ 
&\rightarrow& e^{2m_n} (A_n)^{1/2}[ A_n^{1/2}] = e^{2m_n} A_n
\end{eqnarray*}
since $e^{-m_n} \rightarrow 0$, as $n \rightarrow \infty$. 

Thus $m_n \rightarrow \infty$ implies $A_n \rightarrow 0$. Since the integrand in $A_n$ is positive for finite $\tilde{\sigma}_n$, this implies that $\tilde{\sigma}_n \rightarrow -\infty$
on some open set $W_n \subset M$. But this contradicts that $
\tilde{\sigma}_n \rightarrow \tilde{\sigma} \in C^0(M)$ in $W^{2,2}$, 
lemma (2.3). Thus $m_n$ cannot tend to $\infty$

Note that $m_n$ can still go to $-\infty$, thereby making 
$e^{\sigma_n} \rightarrow u =0$. We will show later that this doesnot happen.

In order to prove (b), set $\psi_n=\Delta_h \tilde{\sigma}_{l_n}$
and observe that, according to  part (a) of Lemma 2.2, the
sequence $\{\psi_n\}$ is bounded in $L^2$. Therefore, passing to a
subsequence, if necessary, there exists $f\in L^2(M)$ such that
\begin{equation*}
\lim_{n\rightarrow\infty}\int_M\psi_n g=\int_M fg
\end{equation*}
for all $g\in L^2(M)$. In particular, considering $g\in C^{\infty}(M)$, this
implies $f =\Delta_h^{distr} \tilde{\sigma}$.

In order to prove (c) we use the following lemma.
\begin{lemma} If a sequence $\{\psi_n\}$ converges to $f\in L^2$ in the weak
topology, then
\begin{equation*}
\lim_{n\rightarrow\infty}\|\psi_n\|\geq\|f\|.
\end{equation*}
Further $lim_{n \rightarrow \infty} \| \psi_n \| = \|f \| $ iff there is
strong convergence.
\end{lemma}
\begin{proof}
The lemma follows from considering the following inequality:
\begin{equation*}
\lim_{n \rightarrow \infty} \int ( \psi_n - f )^2 d \mu \geq 0.
\end{equation*}
\end{proof}
To continue with the proof of the proposition, suppose $\lim_{n\rightarrow
\infty}\|\psi_n\|>\|f\|$;  Using the definition of the functional, we have
\begin{eqnarray*}
S(\sigma_n)
&=& \int_M (K_0 - \frac{1}{2}\Delta_h\tilde{\sigma_n} -Ke^{\tilde{\sigma_n} +m_n})^2 d\mu \\
&=& \frac{1}{4}\|\psi_n\|^2+\|K_0 -Ke^{\tilde{\sigma}_n+ m_n}\|^2 -
\int_M \psi_n(K_0-Ke^{\tilde{\sigma}_n+ m_n})d\mu.
\end{eqnarray*}
From parts (a) and (b) it follows that the sequence $S(\sigma_n)$
converges to $S_0$ and
\begin{eqnarray*}
S_0 &=&\lim_{n\rightarrow\infty} S(\sigma_n)\\
 &=&\lim_{n\rightarrow \infty}\frac{1}{4}\|\psi_n\|^2+\|K_0-Ku\|^2-\int_M f(K_0-Ku)d\mu \\
&>&\frac{1}{4}\|f\|^2+\|K_0-Ku\|^2-\int_M f(K_0-Ku) d\mu \\
&=&\|-\frac{1}{2}f+K_0-Ku\|^2.
\end{eqnarray*}

We will show that this inequality contradicts that $\{\sigma_n\}$ was a
minimizing sequence, i.e.  we can construct a sequence
$\{\tau + m_{l_n} \}$ $ \in$ $  C^{\infty}(M)$
such that $S(\tau + m_{l_n})$ gets as close to $\|-\frac{1}{2}f+K_0-Ku\|^2$
as we like.

Namely, for any $\epsilon>0$ we can construct,
by the density of $C^{\infty}$ in $W^{2,2}$, a function
$\tau\in C^{\infty}(M)$  approximating
$\tilde{\sigma} \in W^{2,2}$  such that $\|\Delta_h\tau-f\|<\epsilon$ and $\|4(v-u)\|<
\epsilon/2$ where $v=\lim_{n\rightarrow \infty}e^{\tau+m_{l_n}}$.
 Since
\begin{equation*}
S_{\tau} = \lim_{n \rightarrow \infty} S(\tau + m_{l_n})=\|-\frac{1}{2}\Delta_h\tau+K_0 -Kv\|^2,
\end{equation*}
we have
\begin{equation*}
|\sqrt{S_{\tau}}-\|-\frac{1}{2}f+K_0-Ku\||\leq \|\frac{1}{2}(f-\Delta_h\tau)-K(v-u)\|\leq\epsilon.
\end{equation*}
Now setting $\delta=\sqrt{S_0}-\|-\frac{1}{2}f+K_0-Ku\|>0$ and choosing
$\epsilon <\delta/2$, and using $\sqrt{S_{\tau}}\leq\|-\frac{1}{2}f+K_0-Ku\|+
\epsilon$ we get, $\sqrt{S_{\tau}}<\sqrt{S_0}-\frac{\delta}{2}$ ---a contradiction,
since $S_0$ was the infimum of the functional.

Thus,
$\lim_{n\rightarrow\infty}\|\Delta_h\tilde{\sigma_n}\|=\|f\|$, so
that, in fact, by lemma 2.5 , the convergence is in the strong
$L^2$ topology. This proves part (c).
\end{proof}

\subsection{Convergence and the non-degeneracy}

\begin{proposition} The minimizing sequence $\{\sigma_{n}\}_{n=1}^{\infty}$
contains a subsequence that converges in $C^0(M)$ to a function $\sigma\in
C^0(M)$, so that the resulting metric $e^{\sigma}h $ is non-degenerate.
\end{proposition}
\begin{proof} Since $\sigma_n=\tilde{\sigma}+m_n$, by proposition $(2.4)$
and lemma $(2.3)$, it is enough to show that the sequence $\{m_n \}$ is
bounded below. Supposing the contrary and passing, if necessary, to a
subsequence, we can assume that
\begin{equation*}
\lim_{n\rightarrow\infty}m_n=-\infty,
\end{equation*}
so that, in notations of proposition $(2.4)$, $u=0$. By proposition
$(2.4),(c)$ we get
\begin{eqnarray*}
S_0=\lim_{n\rightarrow \infty}S(\sigma_{n})=\int_M(K_0-\frac{1}{2}\Delta_h^{distr}
\tilde{\sigma})^2 d\mu.
\end{eqnarray*}
We shall show that this contradicts the fact that $S_0$ is the infimum of the
functional $S$ and that $\{\sigma_n\}$ is a minimizing sequence. First we
have the following lemma.

\begin{lemma} Let $b=K_0-\frac{1}{2}\Delta_h^{distr}\tilde{\sigma}\in L^2(M)$,
where $\tilde{\sigma}_n \rightarrow \tilde{\sigma}$ and $m_n \rightarrow -
\infty$ as $n\rightarrow\infty$. Then
\begin{equation*}
\int_M b\Delta_h\beta d\mu=0
\end{equation*}
for all $\beta\in W^{2,2}(M)$ and  $b \equiv -L$ , where $L$ is a
positive constant.
\end{lemma}

\begin{proof} Consider $G_n(t)=S(\sigma_n+t\beta)-S_0$---a smooth function of
$t$ for a fixed $\beta$. Then by proposition $(2.4), (c)$  we have
\begin{eqnarray*}
G(t) &=& \lim_{n \rightarrow \infty} G_n (t) \\
     &=& \int_M (K_0-\frac{1}{2}\Delta_h^{distr}(\tilde{\sigma} + t\beta))^2
d\mu -\int_M(K_0-\frac{1}{2}\Delta_h^{distr}\tilde{\sigma })^2 d\mu,
\end{eqnarray*}
and $G(t)$ is a smooth function of $t$ for fixed $\beta$. Since $S_0$ is the
infimum of $S$, we have that $G(t) \geq 0$ for all $t$ and $G(0) = 0$.
Therefore it follows that
\begin{equation*}
\frac{dG}{dt}|_{t=0}=0
\end{equation*}
for all $\beta\in W^{2,2}(M)$. Straightforward computation yields
\begin{equation*}
\frac{dG}{dt}|_{t=0}=-\int_{M} b\Delta_h\beta d\mu.
\end{equation*}
Therefore, $b$ satisfies the Laplace equation $\Delta_h b = 0$ in
a distributional sense and from elliptic regularity it follows
that $b$ is smooth. Thus $b$ is harmonic and therefore is a
constant. Finally, by the Gauss-Bonnet theorem, we have $\int_M b
d \mu = 4\pi (1-g) $ and recalling that $g>1$, we conclude that
$b= 4 \pi (1-g) = -L<0$.
\end{proof}

To complete the proof of the proposition, we get a contradiction as follows.
By lemma $(2.7)$ we have that $S_0=\int_{M}(-L)^2 d\mu=L^2$ is the infimum of
the functional. Since $L>0$, and $\{m_n \} \rightarrow -\infty$, we
consider $\tau=\tilde{\sigma}+m_n$ and
choose $n$ large enough so that $-Ke^{\tau}<L/2$. We have
\begin{equation*}
S(u)=\int_M(K_0-\frac{1}{2}\Delta_h \tilde{\sigma}-Ke^{\tau})^2d\mu=\int_M
(-L-Ke^{\tau})^2d\mu.
\end{equation*}
Then, since $ -L+\alpha<-L-Ke^{\tau}<-L/2 $, where $\alpha>0$ is the infimum of
$-Ke^{\tau}$, we have $(-L -Ke^{\tau})^2 <(L-\alpha)^2$ so that
$S(\tau)<L^2$---a contradiction.
\end{proof}

\section{Smoothness and uniqueness}

Here we complete the proof of the main theorem 3.1 by showing that
\begin{proposition} The minimizing function $\sigma\in C^0(M)$ is smooth and
corresponds to the unique K\"{a}hler metric of  negative curvature $K$.
\end{proposition}
\begin{proof}
Let $b=(K_0-\frac{1}{2}\Delta_h^{distr}\sigma-Ke^{\sigma})\in L^2(M)$;
according to the proposition $(2.4), (c)$ and proposition $(2.6)$, $S_0=
\int_M b^2d\mu$. Set $G(t)=S(\sigma+t\beta)-S_0$, where
$\beta\in W^{2,2}(M)$. Repeating arguments in the proof of lemma $(2.7)$,
we conclude that $G(t)$ for fixed $\beta$ is smooth, $G(0)=0$ and $G(t)\geq 0$
for all $t$. Therefore,
\begin{equation*}
\frac{dG}{dt}|_{t=0}=0.
\end{equation*}
A simple calculation yields
\begin{equation*}
\frac{dG}{dt}|_{t=0}=\int_M(-b\Delta_h\beta-2Ke^{\sigma}b\beta) d\mu.
\end{equation*}
Thus $b\in L^2(M)$ satisfies, in a distributional sense, the following equation
\begin{equation} \label{eq}
 -\Delta_h b-2Ke^{\sigma} b=0.
\end{equation}

First, we will show that $b=0$ is the only weak $L^2$ solution to the equation
\eqref{eq}. Indeed, by elliptic regularity $b$ is smooth, so that multiplying
\eqref{eq} by $b$ and integrating over $M$ using the Stokes formula, we get
\begin{equation*}
 \int_M \d b\wedge \ast\d b+\int_Mb^2e^{\sigma}d\mu=0,
\end{equation*}
which implies that $b=0$. Thus we have shown that $S_0=0$.

Second, equation $b=0$ for the minimizing function $\sigma\in C^0(M)$ reads
\begin{equation} \label{boot}
\frac{1}{2}\Delta_h^{distr}\sigma=K_{0}-Ke^{\sigma}\in C^{0}(M).
\end{equation}
Therefore, $\Delta_h^{distr}\sigma$ belongs to $L^p(M)$ so that
$\sigma\in W^{2,p}$ for all $p$. By the Sobolev embedding theorem it follows
that $\sigma$ $\in$ $C^{1,\alpha}(M)$ for some $0<\alpha<1$.Therefore,
the right hand side of the equation \eqref{boot} actually belongs to the space
$C^{1,\alpha}(M)$, and therefore $\sigma\in C^{3,\alpha}(M)$ and so on. This
kind of bootstrapping argument shows that $\sigma$ is smooth~\cite{K}.

The equation $b \equiv 0$ satisfied by $\sigma$ now translates to
$K(\sigma) \equiv K$, where $K(\sigma)$ is the Gaussian curvature of the
metric $e^{\sigma}h dz \otimes d \bar{z}$ , $\sigma \in C^{\infty}(M)$.

The minimizing function $\sigma$ is unique: here is the standard argument,
which goes back to Poincar\'{e}.  Let $\eta$ be another minimizing function,
which is smooth and also satisfies the equation \eqref{boot}
\begin{equation*}
\frac{1}{2}\Delta_h \eta=K_0-Ke^{\eta},
\end{equation*}
so that
\begin{equation*}
\Delta_h (\sigma-\eta)=-2K(e^{\sigma}-e^{\eta}).
\end{equation*}
Multiplying this equation by $\sigma-\eta$ and integrating over $M$ with the
help of Stokes formula, we get

\begin{equation*}
-\int_M \d\xi \wedge \ast\d\xi=\int_M -2K(\sigma - \eta)(e^{\sigma} - e^{\eta})
d\mu,
\end{equation*}
where we set $\xi=\sigma-\eta$. Since $-2K(\sigma - \eta)(e^{\sigma} - e^{\eta})\geq 0$   , we conclude that $\d\xi=0$ and, in fact, $\xi=0$.
\end{proof}

The proof of Theorem 3.1 is complete.

{\bf Acknowledgement}

I would like to thank Professor L. Takhtajan for his invaluable
help in a previous version of this paper. I would also like to
thank Professor D. Geller for his careful reading of the
manuscript. This work was done while I was at SUNY at Stony Brook,
USA.


\begin{thebibliography}{99}
\bibitem{A}T.Aubin: Nonlinear Analysis on manifolds;Monge-Ampere equations.
\bibitem{B} M. S. Berger: Riemannian structures of prescribed Gaussian
Curvature for compact $2$-manifolds J.Differential Geom.(1971)(pg 325-332)
\bibitem{CL} W. Chen, C. Li: A necessary and sufficient condition for the
Nirenberg problem; Comm. Pure. Appl. Math, vol 48 (1995),no. 6, 657-667.
\bibitem{FKra} Farkas and Kra: Riemann Surfaces.
\bibitem{H} N.J. Hitchin: Self-Duality Equations over a Riemann Surface;
Proc. London Math Soc (3) 55 (1987), 59-126.
\bibitem{KW}J. Kazdan, F.W. Warner: Curvature functions for compact
$2$-manifolds; Ann. of Math. (2) 99 (1974), 14-47.
\bibitem{K} J. Kazdan: Applications of PDE to problems in geometry.
\bibitem{P} H. Poincar\'{e}: Les fonctions fuchsiennese l' \'{e}quation
$\Delta u = e^{u}$, J. Math.Pures Appl.(5)4(1898).
\bibitem{U} K. Uhlenbeck: Connections with $L^p$ bounds on
curvature, Comm. Math. Phys. 83 (1982) 31-42.
\bibitem{XY}X.Xu, P.C. Yang: Remarks on Prescribing Gauss curvature;
Trans. AMS, vol 336, no.2, 1993, 831-840.
\end{thebibliography}
\end{document}